\documentclass[12pt]{article}
\usepackage{amsthm}
\usepackage{amsmath}
\usepackage{amsfonts}
\usepackage{cancel}
\usepackage{amssymb}

\usepackage[all, knot]{xy}

\newcommand {\emptycomment}[1]{}

\newcommand{\dM}{\mathrm{d}}

\newcommand{\br}[1]{[\cdot,\cdot]}

\newcommand{\maps}{\colon}

\newcommand{\h}{{H}}
\newcommand{\g}{{Z}}

\newcommand{\calA}{{\cal A}}
\newcommand{\calC}{{\cal C}}
\newcommand{\calD}{{\cal D}}

\newcommand{{\calJ}}{{\cal J}}

\newcommand{\calV}{{\cal V}}

\newcommand{\calZ}{{\cal Z}}

\newcommand{\id}{\mbox{id}}

\newcommand{\Ker}{\operatorname{Ker}}

\newcommand{\half}{\textstyle{\frac{1}{2}}}

\newcommand{\six}{\textstyle{\frac{1}{6}}}

\newcommand{\trl}{\triangleleft}
\newcommand{\trr}{\triangleright}

\newcommand{\pf}{\noindent{\bf Proof.}\ }
\allowdisplaybreaks

\newtheorem{Theorem}{Theorem}[section]
\newtheorem{Proposition}[Theorem]{Proposition}
\newtheorem{Definition}[Theorem]{Definition}

\newtheorem{Example}[Theorem]{Example}

\title{Zinbiel 2-algebras}

\author{Tao Zhang}

\date{}

\begin{document}


\maketitle
 \setcounter{section}{0}

 \vskip0.1cm

{\bf Abstract}\quad
The notions of Zinbiel 2-algebras  and 2-term $Z_{\infty}$-algebras  are introduced and studied.
It is proved that the category of Zinbiel 2-algebras and the category of $2$-term $Z_{\infty}$-algebras are equivalent.
Relationship between Zinbiel 2-algebras and dendriform 2-algebras are found.

\tableofcontents

\section{Introduction}


The notion of $L_{\infty}$-algebras which generalizes Lie algebras first appeared in deformation theory and then in closed string field theory.
The algebraic theory of Lie 2-algebras was studied by Baez and Crans in [1].
It is showed that a Lie 2-algebra can be seen as a categorification of a Lie algebra, where the underling vector space is replaced by 2-vector space and the
Jacobi identity is replaced by a nature transformation which satisfies some coherence law.
They proved that the category of Lie  $2$-algebras and the category of $2$-term $L_{\infty}$-algebras are equivalent.

On the other hand,  the study of Zinbiel algebras was initiated by Loday [2].
Under the Koszul duality the operad of Zinbiel algebras is dual to the operad of Leibniz algebras.
Now a natural question arise, does there exist a categorification of a Zinbiel algebra?
In this paper, we give a positive answer to this question.

The key point is that a Zinbiel 2-algebra is actually a Zinbiel algebra in  the category of 2-vector spaces, where the condition is replaced by some nature transformation.
We also introduced the concept of 2-term $Z_{\infty}$-algebras, and proved that the category of Zinbiel 2-algebras and the category of $2$-term $Z_{\infty}$-algebras are equivalent.

The second aim of this paper is to category the relation between Zinbiel algebras and dendriform algebras. We introduce the notion of Zinbiel 2-algebra, and establish the following commutative diagram:
         $$
\xymatrix{
\mbox{Zinbiel~2-algebras}\ar[r]^{}&\mbox{Dendriform~2-algebras}\ar[r]^{}& \mbox{Ass~2-algebras}\\
\mbox{Zinbiel algebras}\ar[u]_{\mbox{Cat}}\ar[r]&\mbox{Dendriform algebras}\ar[u]_{\mbox{Cat}}\ar[r]& \mbox{Ass algebras.}\ar[u]_{\mbox{Cat}}}
                $$

The main results and organization of this paper are as follows.
In Section 2, we recall some facts about Zinbiel algebras.
In Section 3, we  introduce the notion of Zinbiel 2-algebras and 2-term $Z_{\infty}$-algebras.
It is proved that there is an equivalence between the category of Zinbiel 2-algebras and the category of $2$-term $Z_{\infty}$-algebras.
In Section 4, we investigate some special cases of Zinbiel 2-algebras, such as skeletal and strict ones.
In Section 5, we found that  Zinbiel 2-algebras give rise to dendriform 2-algebras and associative 2-algebras.

\section{Preliminaries}

In this section, we will recall some facts and definitions about Zinbiel algebras.

\begin{Definition}\label{def:colorLie}  A Zinbiel algebra  is a vector space $Z$ together with a multiplication $\cdot: \g \times \g \to \g$ satisfying the following Zinbiel identity:
\begin{align}\label{J1}
(x\cdot y)\cdot z=x\cdot (y\cdot  z)+x\cdot (z\cdot  y)
\end{align}
for all $x,y, z\in \g$.
\end{Definition}

\begin{Example}\label{ex:colorglnv} Let $(Z,\cdot)$ be a Zinbiel algebra with multiplication $\cdot:Z \times Z\to Z$.
Define the  symmetrized product by
$$x y:=\half (x\cdot y+y\cdot  x),$$
for all $x, y\in A$.  Then under the symmetrized product,  $Z$ becomes an associative and commutative algebra.
\end{Example}

A homomorphism  between two Zinbiel algebras $(\g,\cdot)$ and $(\g', \cdot')$ is a linear map $\varphi: \g \to \g'$ such that
$$\varphi(x\cdot y) = \varphi(x)\cdot' \varphi(y)$$ 
for all $x, y\in  \g$.

Let $\g$ be a Zinbiel algebra and  $V$ is a vector space, a {\bf bimodule} of $Z$ over the vector space $V$ is a pair of  bilinear maps $\trr:{Z}\times {V} \to {V}$ and $\trl:{V}\times {Z} \to {V}$ such that the following conditions hold:
\begin{eqnarray}
  &&(x \cdot y) \trr v = x \trr (y \trr v+v\trl y),\\
  && (v \trl x) \trl y =v\trl(x \cdot y+y \cdot x) ,\\
  &&(x\trr v)\trl y  = x\trr (v \trl y+y \trr v),
\end{eqnarray}
for all $ x,y,z \in {Z}$, $u \in V$.

\begin{Example}\label{example1}
 Let ${Z}$ be a 2-dimension vector space with $\{e_1 ,e_2\}$  a base. We defined the multiplications by
\[
e_1\cdot e_1 = e_2,
\quad
e_1\cdot e_2 = e_2\cdot e_1=0,
\quad
e_2\cdot e_2 = 0.
\]
Then $({Z}, \cdot)$ is a 2-dimension Zinbiel algebra.
\end{Example}

Let $Z$ be a Zinbiel algebra and let $V$ be a $Z$-bimodule. Recall that the low dimensional cohomology is  defined as following, see [4] .

The cochain complex is
$$
C^{n}(Z, V)=\operatorname{Hom}_{K}\left(Z^{\otimes n}, V\right)
$$
For $1 \leq i \leq 3$, the coboundary operator
$$
d^{i}: C^{i}(Z, V) \rightarrow C^{i+1}(Z, V)
$$
is defined by
$$
\begin{aligned}
\left(d^{1} \omega\right)(x, y)=& x \trr \omega(y)-\omega(x \cdot y)+\omega(x) \trl y, \\
\left(d^{2} \omega\right)(x, y, z)=& x \trr (\omega(y, z)+\omega(z, y))-\omega(x \cdot y, z)+\omega(x, y \cdot z+z \cdot y)-\omega(x, y) \trl z,\\
\left(d^{3} \omega\right)(x, y, z, t)=& x \trr \{\omega(y, z, t)-\omega(z, t, y)+\omega(z, y, t)-\omega(t, z, y)\} \\
&-\omega(x \cdot y, z, t)+\omega(x, y \cdot z+z \cdot y, t) \\
&-\omega(x, y, z \cdot t+t \cdot z)+\omega(x, y, z) \trl t.
\end{aligned}
$$
It is easy to check directly that $d^{i+1} d^{i}=0$ for $i=1,2$, and so $\left(C^{*}(Z, V), d\right)$, $1 \leq * \leq 4$, is a cochain complex, which is usually abbreviated to $C^{*}(Z, V)$. For $n=2,3$, we get the cohomology group
$\mathbf{H}^{n}({Z}, V)=Z^{n}({Z}, V)/B^{n}({Z}, V)$ where $Z^{n}({Z}, V)$ is the space of cocycles
and $B^{n}({Z}, V)$ is the space of coboundaries.


\section{Zinbiel 2-algebras and 2-term $Z_{\infty}$-algebras}

In this section, we introduced the concept of Zinbiel 2-algebras and 2-term $Z_{\infty}$-algebras.
It is proved that the category of Zinbiel 2-algebras and the category of $2$-term $Z_{\infty}$-algebras are equivalent.

Recall that we denote the category of vector spaces by ${\rm\mathbf{Vect}}$.
A {\bf  2-vector space} is a category in ${\rm\mathbf{Vect}}$.
Thus, a 2-vector space $\calV$ is a category with a vector space of objects $\calV_0$ and a vector space of morphisms $\calV_1 $, such that the
source and target maps $s,t \maps \calV_{1} \rightarrow \calV_{0}$, the
identity-assigning map $i \maps \calV_{0} \rightarrow \calV_{1}$, and the
composition map $\cdot \maps \calV_{1} \times _{\calV_{0}} \calV_{1}
\rightarrow \calV_{1}$ are all  linear maps.  We write a
morphism $f$ from source $x$ to target $y$ by $f \maps x \to y$, i.e. $s(f) = x$ and $t(f) = y$.
We also write $i(x)$ as $1_x$.

It is well known that 2-vector spaces are in one-to-one correspondence with 2-term complexes of vector spaces.
A 2-term complex of vector spaces is a pair of vector space with a differential between them: $\calC_1\stackrel{\dM}
{\longrightarrow}\calC_0$.  The correspondence is given as follows.
Given a 2-vector space $\calV$,  then $\Ker(s)\stackrel{t}{\longrightarrow}\calV_0$ is a 2-term complex.
Conversely, any 2-term complex of vector spaces
$\calV_1\stackrel{\dM}{\longrightarrow}\calV_0$ gives rise to a
2-vector space of which the set of objects is $\calC_0$, the set of
morphisms is $\calC_0\oplus \calC_1$, the source map $s$ is given by
$s(x,h)=x$, and the target map $t$ is given by $t(x,h)=x+dh$,
where $x\in \calV_0,~h\in \calV_1.$ The 2-vector space associated
to the 2-term complex of vector spaces is denoted by
$\calV_1\stackrel{\dM}{\longrightarrow}\calV_0$ by $\calV$:
\begin{equation}\label{eqn:V}
\calV=\begin{array}{c}
\calV_1:=\calC_0\oplus \calC_1\\
\vcenter{\rlap{s }}~\Big\downarrow\Big\downarrow\vcenter{\rlap{t }}\\
\calV_0:=\calC_0.
 \end{array}\end{equation}

\begin{Definition} \label{defnlie2alg}
A {\bf Zinbiel 2-algebra} consists of a  2-vector space $\calZ$ equipped with
a bilinear functor, $\cdot\maps {\calZ} \times {\calZ}\rightarrow {\calZ}$
and a natural isomorphism, the {\bf Zinbielator},
$$J_{x,y,z} \maps (x\cdot y)\cdot z\to x\cdot (y \cdot z)+x\cdot (z \cdot y),$$
such that the following  {\bf Zinbielator identity} is satisfied
\begin{align}\label{Jacobiator}
&\quad \Big(x\cdot J_{y,z,t}+x\cdot J_{z,y,t}+1+1\Big)\Big(J_{x,y\cdot z,t}+ J_{x,z\cdot y,t}\Big)J_{x,y,z}\cdot t\notag\\
&=\Big(1+1+x\cdot J_{z,t,y}+x\cdot J_{t, z,y}\Big)\Big(J_{x,y, z\cdot t}+ J_{x,y, t\cdot z}\Big)J_{x\cdot y,z,t}
\end{align}
In terms of a commutative diagram, we see that it relates two
ways of using the Zinbielator to rewrite the expression $((x\cdot y)\cdot z)\cdot t$:
$$\def\objectstyle{\scriptstyle}
  \def\labelstyle{\scriptstyle}
\xymatrix{
&((x\cdot y)\cdot z)\cdot t
\ar[dr]^{J_{x\cdot y,z,t}}\ar[dl]_{J_{x,y,z}\cdot t}&\\
x\cdot ((y\cdot z)\cdot t)+ (x\cdot (z\cdot y))\cdot t\ar[dd]^{J_{x,y\cdot z,t}+ J_{x,z\cdot y, t}}
&& (x\cdot y)\cdot (z\cdot t)+(x\cdot y)\cdot (t\cdot z)\ar[dd]_{J_{x,y, z\cdot t}+  J_{x,y, t\cdot z}}\\
&&\\
{\begin{aligned}&\scriptstyle x\cdot( (y\cdot z)\cdot t)+ x\cdot( (z\cdot y)\cdot t)\\[-.5em]
&\scriptstyle+ x\cdot(t\cdot (y\cdot z))+ x\cdot (t\cdot (z\cdot y))\end{aligned}}
\ar[dr]_{x\cdot J_{y,z,t}\cdot t+x\cdot J_{z,y,t}+1+1}&&
{\begin{aligned}&\scriptstyle  x\cdot (y\cdot (z\cdot t))+ x\cdot (y\cdot (t\cdot z))\\[-.5em]
&\scriptstyle +  x\cdot ((z\cdot t))\cdot y)+ x\cdot ((t\cdot z)\cdot y)\end{aligned}}
\ar[dl]^{1+1+x\cdot J_{z,t,y}+ x\cdot J_{t,z,y}}\\
&P&}
\\ \\
$$
where $P$ is given by
\begin{eqnarray*}
  P&=&x\cdot (y\cdot (z\cdot t))+x\cdot (y\cdot (t\cdot z))\\
  &&+ x\cdot (z\cdot (y\cdot t))+x\cdot (z\cdot (t\cdot y))\\
   &&+x\cdot (t\cdot (y\cdot z))+x\cdot (t\cdot (z\cdot y)).
\end{eqnarray*}
\end{Definition}

\begin{Definition}
  Given Zinbiel 2-algebras $({\calZ},\cdot)$ and $({\calZ}^\prime,\cdot')$, a Zinbiel 2-algebra homomorphism $F:{\calZ}\longrightarrow  {\calZ}^\prime$ consists of:
  \begin{itemize}
   \item[$\bullet$] a functor $(F_0,F_1)$ from the underlying
   2-vector space of ${\calZ}$ to that of ${\calZ}^\prime$.
 \item[$\bullet$] a natural transformation
 $$
F_2(x,y):F_0(x)\cdot' F_0(y)\longrightarrow F_0(x\cdot y)
 $$
such that the following diagram commutes:
  \end{itemize}
 $$
\footnotesize{ \xymatrix{
 (F_0(x)\cdot' F_0(y))\cdot' F_0(z)\ar[d]_{J_{F_0(x),F_0(y),F_0(z)}}\ar[rr]^{\qquad\qquad\qquad F_2(x,y)
\quad\qquad\quad}&&F_0(x\cdot y)\cdot' F_0(z)\ar[d]^{F_2(x\cdot y, z)}\\
~F_0(x)\cdot'(F_0(y)\cdot'F_0(z))+F_0(x)\cdot'(F_0(z)\cdot'F_0(y)))\ar[d]_{F_2(x,y)+F_2(z,y)}&&F_0((x\cdot y)\cdot z)\ar[d]^{F_1J_{x,y,z}}\\
~F_0(x)\cdot'F_0(y\cdot z)+F_0(x)\cdot'F_0(z\cdot y)\ar[rr]^{\footnotesize{F_2(x\cdot y,z)+F_2(y,x\cdot z)}}&&
\footnotesize{F_0(x\cdot (y\cdot z))+ F_0(x\cdot (z\cdot y)). }}}
$$
\end{Definition}

The identity homomorphism ${\id}_{\calZ}:{\calZ}\longrightarrow {\calZ}$ has the identity functor as its underlying functor, together with an identity natural
transformation as $({\id}_{\calZ})_2$.  Let ${\calZ},~{\calZ}'$ and ${\calZ}''$ be Zinbiel 2-algebras, the composition of a pair of Zinbiel 2-algebra morphisms
$F:{\calZ}\longrightarrow {\calZ}'$ and $G:{\calZ}'\longrightarrow {\calZ}''$, which we denote by $G\cdot F$,
is given by letting the functor $((G\cdot F)_0,(G\cdot F)_1)$ be the usual composition of $(G_0,G_1)$ and $(F_0,F_1)$, and letting $(G\cdot F)_2$ be the following composite:
$$
\footnotesize{ \xymatrix{
(G_0\cdot F_0(x))\cdot'' (G_0\cdot F_0(y))\ar[dd]_{G_2(F_0(x),F_0(y))}\ar[dr]^{(G\cdot F)_2(x,y)}&&\\
&G_0\cdot F_0(x\cdot y).&\\
 G_0(F_0(x)\cdot' F_0(y))\ar[ur]_{G_1(F_2(x,y))}&&}}
$$
It is easy to see that there is a category {\bf Z2Alg} with Zinbiel 2-algebras as objects and Zinbiel 2-algebra homomorphisms as morphisms.

\begin{Definition} \label{2termJordan}
A $2$-term $Z_{\infty}$-algebra $\calV=\calV_0\oplus \calV_1$ is a complex consisting of the following
data:
\begin{itemize}
  \item two vector spaces $\calV_{0}$ and
   $\calV_{1}$ together with a linear map    $l_1=d\maps \calV_{1} \rightarrow \calV_{0}$,

  \item a bilinear map $l_{2}\maps \calV_{i} \times \calV_{j}
   \rightarrow \calV_{i+j},$ where $0 \leq i + j \leq 1$,

  \item a trilinear map $l_{3}\maps \calV_{0} \times \calV_{0} \times
   \calV_{0} \rightarrow \calV_{1}.$
\end{itemize}
These maps satisfy the following conditions: for all $x,y,z,t\in \calV_{0}$ and $h, k \in \calV_{1}$,
  \item[(a)]  $l_2(h,k)=0$,
  \item[(b1)] $\dM (l_2(x,h)) = l_2(x,\dM h)$,
 \item[(b2)] $\dM (l_2(h,x)) = l_2(\dM h,x)$,
  \item[(c)] $l_2(\dM h,k) =  l_2(h,\dM k)$,
  \item[(d )] $\dM (l_{3}(x,y,z))=l_2(x,l_2(y,z))+l_2(x,l_2(z,y))-l_2(l_2(x,y),z)$,
  \item[(e1)] $l_{3}(x,y,\dM h)=l_2(x,l_2(y,h))-l_2(x,l_2(y,h))-l_2(x,l_2(h,y))$,
 \item[(e2)] $l_{3}(x,\dM h,z)=l_2(l_2(x, h), z)-l_2(x,l_2(h, z))-l_2(x,l_2(z, h))$,
  \item[(e3)] $l_{3}(\dM h,y,z)=l_2(l_2(h,y), z)-l_2(h,l_2(y, z))-l_2(h,l_2(z, y))$,
  \item[(f)] $l_2(x, l_3(y, z, t))+l_2(x, l_3(z, y, t)) +l_3(x, l_2(y, z)+l_2(z,y), t) +l_2(l_3(x, y, z), t)$\\
      $-l_2(x,  l_3(z, t, y))-l_2(x, l_3(t, z, y))-l_3(x, y, l_2(z,t)+l_2(t, z))-l_3(x \cdot y, z, t)=0$.
\end{Definition}

\begin{Definition}\label{defi:Lie-2hom}
Let $(\calV;\dM,l_2,l_3)$ and $(\calV';\dM',l'_2,l'_3)$ be two $2$-term $Z_{\infty}$-algebras.
A $Z_\infty$-homomorphisms $f$ from $\calV$ to $\calV'$ consists of
 linear maps $f_0:\calV_0\rightarrow \calV_0',~f_1:\calV_{1}\rightarrow \calV_{1}'$
 and $f_{2}: \calV_{0}\times \calV_0\rightarrow \calV_{1}'$,
such that the following equalities hold for all $ x,y,z\in \calV_{0},
a\in \calV_{1},$
\begin{itemize}
\item [$\rm(i)$] $f_0\dM=\dM'f_1$,
\item[$\rm(ii)$] $f_{0}l_2(x,y)-l'_2(f_{0}(x),f_{0}(y))=\dM'f_{2}(x,y),$
\item[$\rm(iii)$] $f_{1}l_2(x,a)-l'_2f_{0}(x),f_{1}(a))=f_{2}(x,\dM a)$,
\item[$\rm(iv)$] $f_1(l_3(x,y,z))-l_3'(f_0(x),f_0(y),f_0(z))$
  $=f_2(x, l_2(y,z))- f_2(l_2(x,y),z) - f_2(y,l_2(x,z))$ $+ l'_2(f_0(x), f_2(y,z)) - l'_2(f_2(x,y), f_0(z))- l'_2(f_0(y), f_2(x,z)).$
\end{itemize}
 If $f_2=0$, the $Z_\infty$-homomorphism $f$ is called a strict $Z_\infty$-homomorphism.
\end{Definition}

Let $f:\calV\to \calV'$ and $g:\calV'\to \calV''$ be two $Z_\infty$-homomorphisms, then their composition $g\cdot f:\calV\to \calV''$
is a $Z_\infty$-homomorphism defined as $(g\cdot f)_0=g_0\cdot f_0$, $(g\cdot f)_1=g_1\cdot f_1$
and
$$(g\cdot f)_2(x,y)=g_2 (f_0(x), f_0(y))+g_1(f_2(x,y)).$$

The identity $Z_\infty$-homomorphism $1_\calV: \calV\to \calV$ has the identity chain map together with $(1_\calV)_2=0$.
There is a category {\bf 2TermZ$_\infty$} with 2-term $Z_{\infty}$-algebras as objects and $Z_\infty$-homomorphisms as morphisms.

Now we establish the equivalence between the category of Zinbiel 2-algebras and $2$-term $Z_{\infty}$-algebras.

\begin{Theorem} The categories {\bf 2TermZ$_\infty$} and  ${\bf Z2Alg}$ are equivalent.
\end{Theorem}

\pf  We give a sketch of the proof. First, we show how to construct a Zinbiel 2-algebra
from a 2-term $Z_{\infty}$-algebra. 

Let $\calV=(\calV_1\stackrel{\dM}{\longrightarrow}\calV_0,l_2,l_3)$ be a 2-term $Z_{\infty}$-algebra, we introduce a bilinear
functor  on the $2$-vector space ${\calZ}=(\calV_0\oplus \calV_1\rightrightarrows \calV_0)$ given by \eqref{eqn:V},
that is ${Z}$ has vector spaces of objects and morphisms ${\calZ}_0=\calV_0$, ${\calZ}_1=\calV_0\oplus \calV_1$ and a morphisms
$f\maps x \rightarrow y$ in ${\calZ}_{1}$ by
$f=(x, h)$ where homogenous elements $x \in \calV_{0}$ and $h \in \calV_{1}$ have the same degree.  The source, target, and
identity-assigning maps in $L$ are given by
\begin{eqnarray*}
  s(f) &=& s(x, h) = x ,\\
  t(f) &=& t(x, h) = x + dh, \\
  i(x) &=& (x, 0),
\end{eqnarray*}
and we have $t(f) - s(f) = dh$.

Now we define product on ${\calZ}$ by
$$(x,h)\cdot (y,k)=l_2(x,y)+l_2(x,k)+l_2(h,y)+l_2(dh,k).$$
It is straightforward to see that it is a bilinear functor.

Now we define the Zinbielator as following
$$J_{x,y,z}:=((x\cdot y)\cdot z,l_3(x,y,z)).$$
Then by Condition $(g)$, we have $J_{x,y,z}$ is a morphism from source $(x\cdot y)\cdot z$ to target $x\cdot (y\cdot z)+x\cdot (z\cdot )$.
One show that $J_{x,y,z}$ is natural isomorphism and satisfying Zinbielator identity \eqref{Jacobiator} in Definition \ref{defnlie2alg}.
Thus from a 2-term $Z_{\infty}$-algebra, we can obtain a Zinbiel 2-algebra.

For any $Z_\infty$-homomorphism $f=(f_0,f_1,f_2)$ form $\calV$ to
$\calV'$, we construct a Zinbiel 2-algebra morphism $F=T(f)$
from $L=T(\calV)$ to $L'=T(\calV')$ as follows:

Let $F_0=f_0,~F_1=f_0\oplus f_1$, and $F_2$ be given by
$$
F_2(x,y)=(f_0(x)\cdot f_0(y),f_2(x,y)).
$$
Then $F_2(x,y)$ is a bilinear natural isomorphism from $F_0(x)\cdot F_0(y)$ to $F_0(x\cdot y)$, and $F=(F_0,F_1,F_2)$ is a Zinbiel 2-algebra homomorphism from ${\calZ}$ to ${\calZ}'$.
One can also deduce that $T$ preserves the identity $Z_\infty$-homomorphisms and the composition of $Z_\infty$-homomorphisms. Thus, $T$ constructed above is a
functor from {\bf 2TermZ$_\infty$} to {\bf Z2Alg}.

Conversely, given a Zinbiel 2-algebra ${\calZ}$, we define $l_2$ and
$l_3$ on the 2-term complex ${\calZ}_1\supseteq \ker(s)=\calV_1\stackrel{\dM}{\longrightarrow}\calV_0={\calZ}_0$
by
\begin{itemize}
\item $\dM h = t(h)$ for $h \in \calV_1 \subseteq {\calZ}_1$.
\item $l_{2}(x,y) = x\cdot y$ for $x,y \in \calV_0 = {\calZ}_0$.
\item $l_{2}(x,h) = [1_x, h]$
for $x \in \calV_0 = \calZ_0$ and $h \in \calV_1 \subseteq {\calZ}_1$.
\item $l_2(h,k) = 0$ for $h,k \in \calV_1 \subseteq {\calZ}_1$.
\item $l_{3}(x,y,z) =p_1J_{x,y,z}$ for $x,y,z \in \calV_0 = {Z}_0$, where $p_1: {Z}_1=\calV_0\oplus \calV_1\longrightarrow \calV_1$ is the projection.
\end{itemize}
Then one can verify that
$(\calV_1\stackrel{\dM}{\longrightarrow}\calV_0,l_2,l_3)$ is a 2-term $Z_{\infty}$-algebra.

Let $F=(F_0,F_1,F_2):{\calZ}\longrightarrow {\calZ}'$ be a Zinbiel 2-algebra homomorphism, and $S({\calZ})=\calV,~S({\calZ}')=\calV'$. Define
$S(F)=f=(f_0,f_1,f_2)$ as follows. Let $f_0=F_0$, $f_1=F_1|_{V_1=\Ker(s)}$ and define $f_2$ by
$$
f_2(x,y)=F_2(x,y)-i(s(F_2(x,y))).
$$
It is not hard to deduce that $f$ is a 2-term $Z_\infty$-algebra homomorphism.
Furthermore, $S$ also preserves the identity morphisms and the
composition of morphisms. Thus, $S$ is a functor from {\bf Z2Alg} to
{\bf 2TermZ$_\infty$}.

We left it to the reader to show that there are natural isomorphisms
$\alpha:T\cdot S\Longrightarrow 1_{{\bf Z2Alg}}$ and $\beta:S\cdot
T\Longrightarrow 1_{{\bf 2TermZ_\infty}}$.
\qed

\section{Construction of Zinbiel 2-algebras}

In this section, special cases and concrete examples of Zinbiel 2-algebras are given.
This include Zinbiel algebras with  3-cocyles, crossed module of  Zinbiel algebras.

\subsection{Skeletal Zinbiel 2-algebras}

A 2-term $Z_{\infty}$-algebra is called {\bf skeletal} if $d=0$.
In this case,
from conditions $(a)$ and $(g)$, we have $\calV_0$ is a Zinbiel algebra.
Conditions $(b)$ and $(h)$ imply that $\calV_1$ is a representation of $\calV_0$ by the action
defined by $x\trr h := l_2(x, h),  h\trl x := l_2(h, x)$. Now condition $(i)$ can be described in terms of a 3-cocycle condition in the
Zinbiel algebra cohomology of $\calV_0$ with values in $\calV_1$.

\begin{Proposition}
Skeletal 2-term $Z_{\infty}$-algebras are in one-to-one correspondence with
quadruples $(\g, V, l_3)$ where $\g$ is a Zinbiel algebra, $V$ is a bimodule of $\g$ and $l_3$ is a 3-cocycle on $\g$ with values in $V$.
\end{Proposition}

%
%

\subsection{Crossed modules of Zinbiel algebras}
Another kind of  2-term $Z_{\infty}$-algebra is called {\bf strict} if $l_3=0$.
This kind of Zinbiel 2-algebras can be described in terms of crossed modules of Zinbiel algebras.

\begin{Definition} Let $(\g,\cdot_{\g})$ and $(\h,\cdot_{\h})$ be two Zinbiel algebras.
A crossed module of Zinbiel algebras is a homomorphism of Zinbiel algebras
$\varphi: \h\to \g$ together with an action of $Z$ on $H$, denoted by $\trr: Z\times H\to H, (x, h)\mapsto x\trr h$ and $\trl: H\times Z\to H, (h, x)\mapsto  h\trl x$,
such that
$$\varphi(x\trr h) = x\cdot_{\g}\varphi(h),\quad\varphi(h\trl x) = \varphi(h)\cdot_{\g} x,\quad\varphi(h)\trr k = h\cdot_{\h} k= h\trl \varphi(k),$$
for all $h, k\in\h, x\in\g$.
\end{Definition}

\begin{Proposition}
There is an one-to-one correspondence between strict 2-term $Z_{\infty}$-algebras and crossed modules of Zinbiel algebras.
\end{Proposition}

\pf Let $\calV_1\stackrel{\dM}{\longrightarrow} \calV_0$ be a $2$-term $Z_{\infty}$-algebra with $l_3=0$. We construct Zinbiel algebras on $\g=\calV_0$ and $\h=\calV_1$ as follows.
The product on $\g$ and $\h$ are defined by
\begin{eqnarray*}
&&h\cdot_{\h} k:=l_2(\dM h,k)=l_2(h,\dM k),\quad\forall~h,k\in \h=\calV_1;\\
&&x\cdot_{\g} y:=l_2(x,y),\quad\forall~x,y\in \g=\calV_0.
\end{eqnarray*}
By conditions in Definition \ref{2termJordan},
it is easy to see that $\cdot_\g$ satisfies the Zinbiel identity, and we have
\begin{eqnarray*}
 &&(h  \cdot_\h k) \cdot_\h g-h  \cdot_\h (k \cdot_\h g)-h  \cdot_\h (g \cdot_\h k)\\
  &=&l_2(\dM h, k) \cdot_\h g- h\cdot_\h l_2(k, \dM  g)- h\cdot_\h l_2(\dM g, k)\\
  &=&l_2(l_2(\dM h , k),  \dM g))-l_2(l_2(\dM h ,  l_2( k,  \dM g))-l_2(l_2(\dM h ,  l_2(\dM g, k))\\
 &=&0.
 \end{eqnarray*}
Thus $(H,\cdot_\h)$ is a Zinbiel algebra.
Now let $\varphi=\dM$, then we have
$$
\varphi(h\cdot_\h k)=\dM (l_2(dh,k))=l_2(\dM h,\dM k)=\varphi(h)\cdot_\g \varphi(k),
$$
which implies that $\varphi$ is a homomorphism of Zinbiel algebras from $H$ to $Z$.

Now define the action map of $\g\times \h\to \h$ by
$$x\trr h:=l_2(x,h), \quad h\trl x:= l_2(h,x) \in  \h,$$
then it is easy to check that
\begin{eqnarray*}
&&\varphi(x\trr h)=\dM(l_2(x,h)) = l_2(x,\dM h)= x\cdot_{\g} \varphi(h),\\
&&\varphi(h\trl x)=\dM(l_2(h,x)) = l_2(\dM h,x)= \varphi(h)\cdot_{\g} x,\\
&&\varphi(h)\cdot k =l_2(\dM h, k)=h\cdot_{\h} k.
\end{eqnarray*}
Therefore we obtain a crossed module of Zinbiel algebras.

Conversely,  a crossed module of Zinbiel algebras gives rise to a 2-term $Z_{\infty}$-algebra with $\dM=\varphi$,
$\calV_0=\g$ and $\calV_1=\h$, where the brackets are given by
\begin{eqnarray*}
~ l_2(x,y)&:=&x\cdot_{\g} y;\\
~l_2(x,h)&:=&x\trr h,\quad l_2(h,x):=h\trl x;\\
~l_2(h,k)&:=&0.
\end{eqnarray*}
The crossed module conditions give various conditions for $2$-term $Z_{\infty}$-algebras with $l_3=0$.
\qed

\subsection{Crossed module extensions of Zinbiel algebras}

In this section we will study crossed module extensions of Zinbiel algebras and show that they are related to the some elements of cohomology group.

\begin{Definition}
 Let $Z$ be a Zinbiel algebra and   ${M}$ be a bimodule of $Z$. A crossed module extension of  $Z$ by $M$ is a short exact sequence
\begin{equation}\label{diagram:exact}
 \xymatrix{
   0  \ar[r]^{} & {M} \ar[r]^{i} & V \ar[r]^{\partial} & S \ar[r]^{\pi} &  Z  \ar[r]^{} & 0 \\
  }
\end{equation}
such that $\partial:V\rightarrow S$ is a crossed module and $M\cong \Ker\partial$, $T\cong \mathrm{coKer} \partial$.
\end{Definition}

\begin{Definition}
 Two crossed module extensions of Zinbiel algebras
$\partial:V\rightarrow S$  and ${\partial}':{V}'\rightarrow {S}'$  are equivalent,
 if there exists a Zinbiel algebras homomorphism $\alpha:V\to{V}'$ and $\beta: S\to {S}'$  such that the following diagram commutes
 \begin{equation}
\xymatrix{
    0\ar[r]^{} & M  \ar@{=}[d]_{} \ar[r]^{i} & V  \ar[d]_{\alpha} \ar[r]^{\partial } & S  \ar[d]_{\beta}  \ar[r]^{\pi } & Z \ar@{=}[d]_{} \ar[r]^{} & 0  \\
   0\ar[r]^{} & M  \ar[r]^{i'} & {V}'  \ar[r]^{{\partial}'} & {S}'  \ar[r]^{\pi' } & Z  \ar[r]^{} & 0.}
   \end{equation}
The set of equivalent classes of extensions of $Z$ by $M$ is denoted by $\mathbf{CExt}(Z,{M})$.
\end{Definition}

Given a crossed module $\partial: V \rightarrow S$, we consider $Z=\operatorname{coKer}(\partial)$ and $M=\operatorname{Ker}(\partial)$
in the category of  vector spaces. The Zinbiel algebra structure of $S$ induces a Zinbiel algebra structure on $Z$ and the $Z$-bimodule structure on $V$ induces a $Z$-bimodule structure on $M$. Thus we obtained a crossed module extension.


%



The following Theorem \ref{thm:crossed} shows that  crossed module extensions are naturally associated with the elements of third cohomology group.

\begin{Theorem}\label{thm:crossed}
 Let $Z$ be a Zinbiel  algebra and $M$ be a bimodule of $Z$.
Then there is a canonical map:
\begin{equation*}
  \xi:\mathbf{CExt}(Z,M)\rightarrow \mathbf{H^3}(Z,M).
\end{equation*}
\end{Theorem}

\proof We define $\xi:\mathbf{CExt}(Z,M)\rightarrow \mathbf{H^3}(Z,M)$ as follows. Given a crossed module extension of $T$ by ${M}$ and choose linear sections $s: Z \rightarrow S$ such that $\pi s=1$ and $q: \operatorname{Im}(\partial) \rightarrow V, \partial q=1$. For $x, y \in Z$, we have
$$\pi(s(x)\cdot s(y)-s(x\cdot y))=\pi s(x) \cdot \pi  s(y)-\pi s(x\cdot y)=x\cdot y-x\cdot y=0$$
and thus $s(x)\cdot s(y)-s(x\cdot y) \in \operatorname{Im}(\partial) .$ Take $g(x, y)=q(s(x)\cdot s(y)-s(x\cdot y)) \in V$ and define
\begin{eqnarray*}
&&\theta_{\mathcal{E}}(x, y, z)\\
&=&x \trr ({g}(y, z)+{g}(z, y))-{g}(x \cdot y, z)+{g}(x, y \cdot z+z \cdot y)-{g}(x, y) \trl z\\
&=&s(x) \cdot ({g}(y, z)+{g}(z, y))-{g}(x \cdot y, z)+{g}(x, y \cdot z+z \cdot y)-{g}(x, y) \cdot s(z).
 \end{eqnarray*}
Since $\partial$ is a map of crossed modules, it follows that
\begin{eqnarray*}
&&\partial\left(\theta_{\mathcal{E}}(x, y, z)\right)\\
&=& \partial \{x \trr  ({g}(y, z)+{g}(z, y))-\partial{g}(x \cdot y, z)+\partial{g}(x, y \cdot z+z \cdot y)-\partial {g}(x, y) \trl z\}\\
&=& s(x) \cdot(s(y)\cdot s(z)-s(y\cdot z)+ (s(z)\cdot s(y)-s(z\cdot y))\\
&&-(s(x \cdot y)\cdot s(z)+s((x \cdot y)\cdot z)\\
&&+s(x)\cdot s(y \cdot z+z \cdot y)-s(x\cdot (y \cdot z+z \cdot y))\\
&& -(s(x)\cdot s(y))\cdot s(z)+s(x\cdot y) \cdot s(z)\\
&=&0.
 \end{eqnarray*}
Thus  $\theta_{\mathcal{E}}(x, y, z) \in M=\operatorname{ker}(\partial)$ and the map $\theta_{\mathcal{E}}: Z^{\otimes 3} \rightarrow M$ is well defined.
One  can  checks that $d^3\left(\theta_{\mathcal{E}}\right)=0$.
 Thus we have defined a linear map $\theta_{\mathcal{E}}: Z^{\otimes 3} \rightarrow M$
which is a cocycle with respect to the coboundary operator $d^3$.
We define $\xi: \mathbf{CExt}(Z, M) \rightarrow \mathbf{H}^{3}(Z, M)$ by taking $\xi(\mathcal{E})$
to be the class of $\theta_{\mathcal{E}}$ in $\mathbf{H}^{3}(Z, M)$.

Next, we have to check that the class of $\theta_{\mathcal{E}}$ in $\mathbf{H}^{3}(Z, M)$ does not depend on the sections $s$ and $q$.
We show first that the class of $\theta_{\mathcal{E}}$ does not depend on the section $s$. Suppose $\bar{s}: Z \rightarrow S$ is another section of $\pi$ and let $\bar{\theta}_{\mathcal{E}}$ be the map defined using $\bar{s}$ instead of $s$. Since $s$ and $\bar{s}$ are both sections of $\pi$ there exists a linear map $h: Z \rightarrow V$ with $s-\bar{s}=\partial h .$ 
Observe that
\begin{align*}
&s(x)\cdot g(y,z)- \overline{s}(x) \cdot \overline{g} (y,z)\\
=~& (s(x) - \overline{s}(x))\cdot g (y,z) +  \overline{s}(x)\cdot  (g - \overline{g})(y,z) \\
=~& \partial h (x) \cdot q (s(y)\cdot s(z) - s(y\cdot z))  +  \overline{s}(x)\cdot (g - \overline{g})(y,z) \\
=~&  h(x)\cdot (s(y)\cdot s(z) - s(y\cdot z) )+  \overline{s}(x)\cdot (g - \overline{g})(y,z).
\end{align*}
We obtain
$$
\begin{aligned}
&\left(\theta_{\mathcal{E}}-\bar{\theta}_{\mathcal{E}}\right)(x, y, z)\\
=& h(x)\cdot (s(y) \cdot s(z)-s(y\cdot z))+\bar{s}(x)\cdot (g-\bar{g})(y, z) \\
&+ h(x)\cdot (s(z) \cdot s(y)-s(z\cdot y))+\bar{s}(x)\cdot (g-\bar{g})(z, y) \\
&-(g-\bar{g})(x\cdot y, z)+(g-\bar{g})(x, y\cdot z+z \cdot y)\\
& -(s(x)\cdot s(y)-s(x \cdot y))\cdot h(z)  -(g-\bar{g})(x, y) \cdot\bar{s}(z)\\
\end{aligned}
$$
where $g(x, y)=q(s(x)\cdot s(y)-s(x\cdot y))$ and $\bar{g}(x, y)=q(\bar{s}(x) \cdot\bar{s}(y)-\bar{s}(x\cdot y))$. We define a
$\operatorname{map} b: Z^{\otimes 2} \rightarrow V$ as follows.
$$
b(x, y)=s(x)\cdot h(y)-h(x\cdot y)+h(x)\cdot s(y)-h(x)\cdot \partial h(y)
$$
Then a easy calculation shows that $\partial b=\partial(g-\bar{g})$, thus $(g-\bar{g}-b)$ is a map from $Z^{\otimes 2}$ to $M$. 
Now we can replace $(g-\bar{g})$ by $b$ without changing the equality in $\mathbf{H}^{3}(Z, M)$ since the difference is the coboundary $d^2(g-\bar{g}-b)$. 
After replacing $(g-\bar{g})$ by $b$, since $\partial: V \rightarrow A$ is a crossed module, we obtain the equality $\left(\theta_{\mathcal{E}}-\bar{\theta}_{\mathcal{E}}\right)(x, y, z)=0$ in $\mathbf{H}^{3}(Z, M)$.
This proves that the class of $\theta_{\mathcal{E}}$ does not depend on the section $s$.

At last, consider two equivalent crossed module extensions
 \begin{equation}
\xymatrix{
    0\ar[r]^{} & M  \ar@{=}[d]_{} \ar[r]^{i} & V  \ar[d]_{\alpha } \ar[r]^{\partial } & S  \ar[d]_{\beta}  \ar[r]^{\pi } & Z \ar@{=}[d]_{} \ar[r]^{} & 0  \\
   0\ar[r]^{} & M  \ar[r]^{i'} & {V}'  \ar[r]^{{\partial}'} & {S}'  \ar[r]^{\pi' } & Z  \ar[r]^{} & 0.}
   \end{equation}
Let $s: Z \rightarrow A$ and $q: \operatorname{Im}(\partial) \rightarrow V$ be sections of $\pi$ and $\partial$ and let $s^{\prime}: Z \rightarrow A^{\prime}$ and $q^{\prime}: \operatorname{Im}\left(\partial^{\prime}\right) \rightarrow V^{\prime}$ be sections of $\pi^{\prime}$ and $\partial^{\prime}$. Then
\begin{eqnarray*}
&&\theta_{\mathcal{E}}(x, y, z)\\
&=&s(x) \cdot ({g}(y, z)+{g}(z, y))-{g}(x \cdot y, z)+{g}(x, y \cdot z+z \cdot y)-{g}(x, y) \cdot s(z)\\
&=&s(x)\cdot q (s(y) \cdot s(z)-s(y\cdot z))+ s(x)\cdot q(s(z) \cdot s(y)-s(z\cdot y)) \\
&&-q(s(x\cdot y)\cdot s(z)-s((x\cdot y)\cdot z)) \\
&&+q(s(x)\cdot s(y\cdot z+z \cdot y)-s(x\cdot (y\cdot z+z \cdot y))\\
&&-q(s(x)\cdot s(y)-s(x \cdot y))\cdot s(z).
 \end{eqnarray*}
Since $\pi^{\prime} \beta s=1$ then $\beta s$ is another section for $\pi^{\prime}$ and therefore we can now replace $s^{\prime}$ by $\beta s$ and we obtain the following equality in
$\mathbf{H}^{3}(Z, M)$.
$$
\begin{aligned}
&\left(\theta_{\mathcal{E}}-\theta_{\mathcal{E}^{\prime}}\right)(x, y, z)\\
=& \beta s(x)\left(\left(\alpha q-q^{\prime} \beta\right)(s(y)\cdot s(z)-s(y\cdot z))\right)\\
&+ \beta s(x)\left(\left(\alpha q-q^{\prime} \beta\right)(s(z)\cdot s(y)-s(z\cdot y))\right)\\
&-\left(\alpha q-q^{\prime} \beta\right)(s(x\cdot y) \cdot s(z) - s((x\cdot y)\cdot z) \\
&+\left(\alpha q-q^{\prime} \beta\right)(s(x)\cdot s(y \cdot z+z \cdot y)-s(x\cdot(y \cdot z+z \cdot y))\\
&-\left(\alpha q-q^{\prime} \beta\right)(s(x) \cdot s(y)-s(x\cdot y))\cdot \beta s(z)
\end{aligned}
$$
Thus $\left(\theta_{\mathcal{E}}-\theta_{\mathcal{E}^{\prime}}\right)(x, y, z)=d^2 \phi(x, y, z)$, where $\phi : \otimes^2 Z \rightarrow M$ is defined by
\begin{eqnarray*}
\phi (x,y) &=&(\alpha q - q' \beta)(s(x)\cdot s(y)-s (x\cdot y)).
\end{eqnarray*}
 This proves that
$\theta_{\mathcal{E}}=\theta_{\mathcal{E}^{\prime}}$ in $\mathbf{H}^{3}(Z, M)$ and that the class of $\theta_{\mathcal{E}}$ in $\mathbf{H}^{3}(Z, M)$ does not depend on the sections $s$ and $q .$ Therefore the map $\xi:\mathbf{CExt}(Z,M)\rightarrow \mathbf{H^3}(Z,M)$ is well defined.
\qed

\section{Relation with dendriform 2-algebras}
\begin{Definition} \label{ass2algebra}
A $2$-term $A_{\infty}$-algebra $\calA=\calA_0\oplus \calA_1$ is a complex consisting of the following
data:
\begin{itemize}
  \item two vector spaces $\calA_{0}$ and
   $\calA_{1}$ together with a linear map    ${m}_1=\dM\maps \calA_{1} \rightarrow \calA_{0}$,
  \item a bilinear map ${m}_{2}\maps \calA_{i} \times \calA_{j}
   \rightarrow \calA_{i+j},$ where $0 \leq i + j \leq 1$,
  \item a trilinear map ${m}_{3}\maps \calA_{0} \times \calA_{0} \times
   \calA_{0} \rightarrow \calA_{1}.$
\end{itemize}
These maps satisfy the following conditions: for all $x,y,z,t\in \calA_{0}$ and $h, k \in \calA_{1}$,
\begin{itemize}
  \item[(a)]  ${m}_2(h,k)=0$,
  \item[(b1)] $\dM ({m}_2(x,h)) = {m}_2(x,\dM h)$,
 \item[(b2)] $\dM ({m}_2(h,x)) = {m}_2(\dM h,x)$,
  \item[(c)] ${m}_2(\dM h,k) =  {m}_2(h,\dM k)$,
  \item[(d)] $\dM ({m}_{3}(x,y,z))={m}_2(x,{m}_2(y,z))-{m}_2({m}_2(x,y),z)$,
  \item[(e1)] ${m}_{3}(x,y,\dM h)={m}_2(x,{m}_2(y,h))-{m}_2(x,{m}_2(y,h))$,
 \item[(e2)] ${m}_{3}(x,\dM h,z)={m}_2({m}_2(x,  h), z)-{m}_2(x,{m}_2(h, z))$,
  \item[(e2)] ${m}_{3}(\dM h,y,z)={m}_2({m}_2(h,y), z)-{m}_2(h,{m}_2(y, z))$,
  \item[(f)] ${m}_2(x,{m}_3(y,z,t))-{m}_3({m}_2(x,y),z,t)+{m}_3(x,{m}_2(y,z),t)$~$-{m}_3(y,z,{m}_2(x,t))$\\
     $+{m}_2({m}_3(x, y,z),t)=0$.
 \end{itemize}
\end{Definition}

A $2$-term $A_{\infty}$-algebra is called  2-term $C_{\infty}$-algebras if ${m}_{2}$ and ${m}_{3}$ are commutative maps.

\begin{Definition}
A $2$-term $\mathrm{Dend}_\infty$-algebra $\calD=\calD_0\oplus \calD_1$ consists of the following data
\begin{itemize}
\item a chain complex $\mu_1=\dM: \calD_1 \xrightarrow{\dM} \calD_0,$
\item bilinear maps $\mu_2 : \mathbb{K} [C_2] \otimes (\calD_i \otimes \calD_j) \rightarrow \calD_{i+j},$
\item a trilinear maps $\mu_3 : \mathbb{K}[C_3] \otimes (\calD_0 \otimes \calD_0 \otimes \calD_0) \rightarrow \calD_1$
\end{itemize}
such that the following identities are hold:
\begin{itemize}
\item[(i)] $\mu_2 ([r]; m,n) = 0$,
\item[(ii)] $\dM \mu_2 ([r]; a,m) = \mu_2 ([r]; a, \dM m)$,
\item[(iii)] $\dM \mu_2 ([r]; m,a) = \mu_2 ([r]; \dM m, a)$,
\item[(iv)] $\mu_2 ([r]; \dM m, n) = \mu_2 ([r]; m, \dM n)$,
\item[(v)]$(\mu_2 \circ_1 \mu_2 - \mu_2 \circ_2 \mu_2 )([s]; a, b, c) =  d\mu_3 ([s];~a, b, c),$
\item[(vi1)] $(\mu_2 \circ_1 \mu_2 - \mu_2 \circ_2 \mu_2 )([s]; a, b, m) =  \mu_3 ([s];~a, b, \dM m),$
\item[(vi2)] $(\mu_2 \circ_1 \mu_2 - \mu_2 \circ_2 \mu_2 )([s]; a,  m, c) =  \mu_3 ([s];~a, \dM m, c),$
\item[(vi3)]  $(\mu_2 \circ_1 \mu_2 - \mu_2 \circ_2 \mu_2 )([s]; m,  b, c) =  \mu_3 ([s];~ \dM m, b, c),$
\item[(vii)] $(\mu_3 \circ_1 \mu_2 - \mu_3 \circ_2 \mu_2 + \mu_3 \circ_3 \mu_2) ([t]; a, b, c, e) =~ (\mu_2 \circ_1 \mu_3 + \mu_2 \circ_2 \mu_3) ([t]; a, b, c, e), $
	\end{itemize}
for all $a, b, c, e \in \calA_0; ~m, n \in \calA_1;~ [r] \in C_2;~ [s] \in C_3$ and $[t] \in C_4.$
\end{Definition}

In the above definition, we define two maps $\mu_2 : \mathbb{K}[C_2] \otimes (A_i \otimes A_j) \rightarrow A_{i+j}$  by
$$ \mu_2 ([r]; a, m) = \begin{cases} a \prec m ~~~ &\mathrm{ if } ~~ [r]=[1] \\
a \succ m ~~~ &\mathrm{ if } ~~ [r] = [2] \end{cases}$$
for $a\in A_i, m\in A_j$. Other maps related to  $[s] \in C_3$ and $[t] \in C_4$ are defined similarly. For details, see [3].


\begin{Proposition}
From a $2$-term $Z_{\infty}$-algebra $\calV=\calV_0\oplus \calV_1$, we can obtain a 2-term $C_{\infty}$-algebras by
\begin{eqnarray*}
m_{1}&=& l_1,\\
m_{2}(x,y)&=&\half(l_{2}(x,y)+l_{2}(y,x)), \\
m_{2}(x,h)&=& \half(l_{2}(x,h)+l_{2}(h,x)), \\
m_{3}(x,y,z) &=& \six(l_{3}(x,y,z)+l_{3}(y,z,x)\\
&&+l_{3}(z,x,y)+(l_{3}(y,x,z)+l_{3}(z,y,x)+l_{3}(x,z,y)).
\end{eqnarray*}
This  2-term $C_{\infty}$-algebras is called a symmetrization of the Zinbiel 2-algebra $\calV$.
\end{Proposition}

\begin{Proposition}
From a $2$-term $Z_{\infty}$-algebra $\calV=\calV_0\oplus \calV_1$, we can obtain a 2-term $\mathrm{Dend}_{\infty}$-algebras by
\begin{eqnarray*}
\mu_{1}&=& l_1,\\
\mu_{2}([1],x,y)&=&l_{2}(x,y),\quad \mu_{2}([2],x,y)=l_{2}(y,x), \\
\mu_{2}([1],x,h)&=&l_{2}(x,h),\quad \mu_{2}([2],x,h)= l_{2}(h,x), \\
\mu_{3}([1],x,y,z) &=&l_{3}(x,y,z),\quad\mu_{3}([2],x,y,z) =l_{3}(y,z,x)\\
\mu_{3}([3],x,y,z) &=&l_{3}(z,x,y).
\end{eqnarray*}
\end{Proposition}

Let $\calA$ be an associate 2-algebra. A Rota-Baxter operator on it is given by linear maps $R_0 : \calA_0\to \calA_0$ and $R_1 : \calA_1\to \calA_1$ satisfying:
\begin{eqnarray*}
m_2(R_0(x_1),R_0(x_2))&=&R_0(m_2(x_1,R_0(x_2))+m_2(R_0(x_1), x_2)),\\
m_2(R_1(x_1),R_1(h))&=&R_0(m_2(x_1,R_1(h))+m_2(R_0(x_1), h)),\\
m_2(R_1(h),R_0(x_1))&=&R_0(m_2(h,R_0(x_1))+m_2(R_1(h), x_1)),\\
m_3(R_0(x_1),R_0(x_2),R_0(x_3))&=&R_0(m_3(x_1,R_0(x_2),R_0(x_3))\\
&&+m_3(R_0(x_1),x_2,R_0(x_3))+m_3(R_0(x_1),R_0(x_2),x_3)).
\end{eqnarray*}
\begin{Proposition}
Let $\calA$ be an 2-term $C_{\infty}$-algebras and $(R_0,R_1)$ be a Rota-Baxter operator on it. Then  $\calA$ is a   2-term $Z_\infty$-algebra, where
\begin{eqnarray*}
l_2(x_1,x_2) &=&m_2(x_1,R_0(x_2))+m_2(R_0(x_1), x_2),\\
l_2(x_1,h) &=&m_2(x_1,R_0(h))+m_2(R_0(x_1), h),\\
l_2(h,x_1) &=&m_2(h,R_0(x_1))+m_2(R_1(h), x_1),\\
l_3(x_1,x_2,x_3)&=&m_3(x_1,R_0(x_2),R_0(x_3))\\
&&+m_3(R_0(x_1),x_2,R_0(x_3))+m_3(R_0(x_1),R_0(x_2),x_3).
\end{eqnarray*}
\end{Proposition}

\subsection*{Acknowledgements}
Caution! This is a primary edition, something should be modified in the future.

\section*{References}

\qquad[1] J. C. Baez and A. S. Crans, \emph{Higher-Dimensional Algebra VI: Lie 2-Algebras}, Theory Appl. Categ.  12(2004), 492--528.\\

[2] J.-L. Loday, Cup product for Leibniz cohomology and dual Leibniz algebras, Math. Scand. 77 (1995), 189--196.\\

[3] A. Das, Cohomology and deformations of dendriform algebras, and Dend${}_\infty$-algebras, arXiv:1903.11802.\\

[4] D. Yau, Deformation of dual Leibniz algebra morphisms, Comm. Algebra, 35(2007), 1369--1378.

\vskip 20pt
\footnotesize{
\noindent College of Mathematics and Information Science,\\
Henan Normal University, Xinxiang 453007, P. R. China;\\
 E-mail address:\texttt{{ zhangtao@htu.edu.cn}}
}

\end{document}